\documentclass[11pt,reqno]{amsart}
\usepackage{amsmath,amssymb,amsthm,mathrsfs,dsfont}
\usepackage{CJK}
\usepackage{amsmath}
\usepackage{dsfont}
\usepackage{mathrsfs}
\usepackage{amsmath,amssymb}
\usepackage{amsfonts}
\usepackage{hyperref}
\usepackage{amsthm}
\usepackage{graphicx}
\usepackage{subfigure}
\usepackage{xcolor}
\usepackage{overpic}

\newtheorem{theorem}{Theorem}[section]
\newtheorem{lemma}{Lemma}[section]
\newtheorem{definition}{Definition}[section]
\newtheorem{proposition}{Proposition}[section]
\newtheorem{corollary}{Corollary}[section]\numberwithin{equation}{section}
\newtheorem{remark}{Remark}[section]

\newfont{\bb}{msbm10 at 12pt}

\newcommand{\bal}{\begin{aligned}}      \newcommand{\eal}{\end{aligned}}
\newcommand{\ba}{\begin{array}}      \newcommand{\ea}{\end{array}}
\newcommand{\bc}{\begin{center}}     \newcommand{\ec}{\end{center}}
\newcommand{\be}{\begin{enumerate}}  \newcommand{\ee}{\end{enumerate}}
\newcommand{\beq}{\begin{eqnarray}}  \newcommand{\eeq}{\end{eqnarray}}
\newcommand{\beQ}{\begin{eqnarray*}} \newcommand{\eeQ}{\end{eqnarray*}}
\newcommand{\bi}{\begin{itemize}}    \newcommand{\ei}{\end{itemize}}
\newcommand{\bt}{\begin{tabular}}    \newcommand{\et}{\end{tabular}}
\newcommand{\bdm}{\begin{displaymath}} \newcommand{\edm}{\end{displaymath}}

\def\qed{\hfill{Q.E.D.}\smallskip}

\begin{document}

\title{\bf Combinatorial Ricci flow on cusped 3-manifolds}
\author{Xu Xu}
\maketitle

\begin{abstract}
Combinatorial Ricci flow on a cusped $3$-manifold is an analogue of Chow-Luo's combinatorial Ricci flow on surfaces
and Luo's combinatorial Ricci flow on compact $3$-manifolds with boundary for finding complete hyperbolic metrics on cusped
$3$-manifolds.
Dual to Casson and Rivin's program of maximizing the volume of angle structures,
combinatorial Ricci flow finds the complete hyperbolic metric on a cusped $3$-manifold by
minimizing the co-volume of decorated hyperbolic polyhedral metrics.
The combinatorial Ricci flow may develop singularities.
We overcome this difficulty by extending the flow through the potential singularities using Luo-Yang's extension.
It is shown that the existence of a complete hyperbolic metric on a cusped $3$-manifold
is equivalent to the convergence of the extended combinatorial Ricci flow, which gives a new
characterization of existence of
a complete hyperbolic metric on a cusped $3$-manifold  dual to Casson and Rivin's program.
The extended combinatorial Ricci flow also provides an effective algorithm for finding complete hyperbolic metrics on cusped $3$-manifolds.
\end{abstract}


\section{Introduction}
If $M$ is the interior of a compact oriented $3$-manifold with boundary consisting of tori,
$M$ is called a cusped 3-manifold, the complete hyperbolic metric on which is proved to be unique up to isometry \cite{M, P}.
By subdividing $M$ into ideal tetrahedra and giving them shapes with hyperbolic ideal tetrahedra,
Thurston \cite{T1} wrote down a system of gluing equations (now named after him) for the complex parameters of the
hyperbolic ideal tetrahedra, the solution of which corresponds to a complete hyperbolic metric on the ideally triangulated manifold $M$.

However, it is difficult to solve Thurston's gluing equation directly in practice.
Motivated by \cite{DV},
Casson and Rivin proposed to solve Thurston's gluing equation by maximizing the volume function on the space of angle structures,
which is a convex polytope of the dihedral angles of the tetrahedra in the triangulation.
The readers can also refer to \cite{Ch,FG,Lac,L3,L5,L6,R2,R4} and others for topics related to Casson and Rivin's program.
The approach of combinatorial Ricci flow on cusped $3$-manifolds introduced by Yang \cite{Y}
is dual to Casson and Rivin's program, in which
one works on decorated hyperbolic polyhedral metrics
to ensure that there is no shearing around edges and finds the
complete hyperbolic metrics on cusped 3-manifolds by adjusting the cone angles through adjusting the edge lengths
according to the combinatorial Ricci curvature along the edges.

Suppose $M$ is an ideally triangulated cusped 3-manifold with a triangulation $\mathcal{T}$ and $E$ is the edge set in $\mathcal{T}$.
A decorated hyperbolic polyhedral metric on $(M, \mathcal{T})$
is defined to be the hyperbolic conic metric obtained by isometrically gluing decorated hyperbolic ideal
tetrahedra along the faces such that the decorations are preserved.
The decorated hyperbolic polyhedral metrics on $(M, \mathcal{T})$ are parameterized by the edge lengths.
The combinatorial Ricci curvature $K_i$ along an edge $i\in E$ is defined to be
$2\pi$ less the cone angle along the edge $i\in E$.

Motivated by Chow-Luo's combinatorial Ricci flow on surfaces \cite{CL}
and Luo's combinatorial Ricci flow on compact 3-manifolds with boundary \cite{L2},
Yang \cite{Y} introduced the following combinatorial Ricci flow on cusped 3-manifolds
\begin{equation}\label{combinatorial Ricci flow}
\begin{aligned}
\frac{dl_{i}}{dt}=K_{i},
\end{aligned}
\end{equation}
where $l: E\rightarrow \mathbb{R}$ is the edge length function for the decorated hyperbolic polyhedral metric.
The combinatorial Ricci flow (\ref{combinatorial Ricci flow}) may develop singularities,
which correspond to flat hyperbolic ideal tetrahedra \cite{Y}.
Under an additional condition that the combinatorial Ricci flow (\ref{combinatorial Ricci flow})
does not develop singularities along the flow in finite and infinite time,
Yang \cite{Y} proved the convergence of the combinatorial Ricci flow (\ref{combinatorial Ricci flow}).
This condition is too strong to be satisfied in practice.
In this paper, we shall remove this condition by extending the combinatorial Ricci flow (\ref{combinatorial Ricci flow})
through the potential singularities of the flow.
This gives rise to a new characterization of existence of a
complete hyperbolic metric on a cusped 3-manifold using the extended combinatorial
Ricci flow, which is dual to Casson-Rivin's program.
The extended combinatorial Ricci flow also provides
an effective algorithm for finding the complete hyperbolic metric on a cusped $3$-manifold.

\begin{theorem}\label{main theorem introduction}
Suppose $(M, \mathcal{\mathcal{T}})$ is an ideally triangulated cusped $3$-manifold which supports a complete hyperbolic metric.
Then a decorated hyperbolic polyhedral metric $l^*$ corresponds to a complete hyperbolic metric on $(M, \mathcal{T})$ if and only if the extended
combinatorial Ricci flow converges exponentially fast to $l^*$ for any initial decorated hyperbolic polyhedral metric on $(M, \mathcal{\mathcal{T}})$.
\end{theorem}

The decorated hyperbolic metric $l^*$ in Theorem \ref{main theorem introduction} has zero combinatorial Ricci curvature
and was proved to be unique up to change of decorations by Luo-Yang \cite{LY}
using nonnegative angle structures and Fenchel dual of the volume function.
We will provide a short and direct proof of Luo-Yang's rigidity without involving
angle structures and Fenchel dual.
A result paralleling to Theorem \ref{main theorem introduction}
has been obtained by the author \cite{X3} for compact 3-manifolds with boundary, which
confirms a conjecture of Luo \cite{L2} on the global convergence of combinatorial Ricci flow
on compact 3-manifolds with boundary.
There have been many important work on combinatorial curvature flows on surfaces. See \cite{CL,Ge,GH0,GX2,GGLSW,GLSW,L1,ZX} and others.

The paper is organized as follows.
In Section \ref{Section 2}, we present preliminary results on decorated hyperbolic ideal tetrahedra and decorated
hyperbolic polyhedral metrics on cusped 3-manifolds.
In Section \ref{Section 3}, we study the basic properties of the combinatorial Ricci flow (\ref{combinatorial Ricci flow})
for decorated hyperbolic polyhedral metrics on cusped 3-manifolds.
In Section \ref{Section 4},
we extend the definition of dihedral angles in a decorated hyperbolic ideal tetrahedron and the definition of
combinatorial Ricci curvature on cusped 3-manifolds.
In Section \ref{Section 5}, we extend the
combinatorial Ricci flow (\ref{combinatorial Ricci flow}) through the potential singularities
and prove a generalization of Theorem \ref{main theorem introduction}.
%In Section \ref{Section 5}, we discuss the prescribing combinatorial Ricci curvature problem.

\section{Preliminary on decorated hyperbolic ideal tetrahedra and cusped 3-manifolds}\label{Section 2}
In this section, we give some preliminaries on ideal triangulations of cusped $3$-manifolds,
decorated hyperbolic ideal tetrahedra and decorated hyperbolic polyhedral metrics on cusped $3$-manifolds.
For more details, please refer to \cite{BB,FG,L5,LY, R1,R3}.
\subsection{Triangulations of cusped 3-manifolds}
Suppose $(M^*, \mathcal{T}^*)$ is a closed $3$-manifold with a triangulation $\mathcal{T}^*=(V, E, F, T)$,
where $V, E, F, T$ are the sets of vertices, edges, faces and tetrahedra respectively.
If we remove the vertices from $\mathcal{T}$, then we get an ideal triangulation $\mathcal{T}$
of the manifold $M=M^*\setminus V$, which is composed of ideal tetrahedra.
If for each vertex $p\in V$, the link of $p$ in $M^*$ is a torus, then
$(M, \mathcal{T})$ is called an ideally triangulated cusped $3$-manifolds.
We still call $E, F, T$ as the sets of vertices, edges and tetrahedra in $\mathcal{T}$ respectively.
Each vertex $p\in V$ corresponds to a cusp of the ideally triangulated cusped $3$-manifold $(M, \mathcal{T})$.
For an ideally triangulated cusped $3$-manifold $(M, \mathcal{T})$,
the number of edges is equal to the number of tetrahedra in the triangulation.
For simplicity, we will denote $N$ as the number of the edges and $s$ as the number of cusps in $\mathcal{T}$ in the following.

\subsection{Decorated hyperbolic ideal tetrahedron}
A hyperbolic ideal tetrahedron is a hyperbolic tetrahedra in $\mathbb{H}^3$
with vertices at infinity. The vertices at infinity are different and not part of the
hyperbolic ideal tetrahedron.
A hyperbolic ideal tetrahedron could be taken as the convex hull of four distinct points $v_1, v_2, v_3, v_4$ in $\mathbb{H}^3$
such that $v_1, v_2, v_3, v_4$ are not in a round circle.

Similar to the decorated ideal triangles in \cite{Penner}, one can also introduce the
decoration to parameterize ideal tetrahedron by edge length.
Suppose $\tau$ is a hyperbolic ideal tetrahedron generalized by $v_1, v_2, v_3, v_4\in \partial \mathbb{H}^3$.
$H_1, H_2, H_3, H_4$ are four horospheres attached to the ideal vertices $v_1, v_2, v_3, v_4$ respectively.
Then $\sigma=(\tau, (H_1, H_2, H_3, H_4))$ is called a decorated hyperbolic ideal tetrahedron and $(H_1, H_2, H_3, H_4)$
is called a decoration.
For the edge $e=v_iv_j$ of the hyperbolic ideal tetrahedron $\tau$, if
$H_i\cap H_j=\emptyset$, we define the edge length $l_{ij}$ of $v_iv_j$ to be the hyperbolic distant of $H_i\cap e$ and $H_j\cap e$.
If $H_i\cap H_j\neq\emptyset$, then $-l_{ij}$ is the hyperbolic distant of $H_i\cap e$ and $H_j\cap e$.
Using Penner's cosine law, the length of the arc between $H_i\cap v_iv_j$ and $H_i\cap v_iv_k$ in $H_i\cap \triangle ijk$
is given by $e^{(l_{jk}-l_{ij}-l_{ik})/2}$.
Note that $H_i\cap \tau$ is a Euclidean triangle.
Using Penner's cosine law for decorated hyperbolic ideal triangles, we have the following characterization of nondegenerate decorated hyperbolic ideal
tetrahedra in terms of the edge lengths.

\begin{lemma}[\cite{BPS} Lemma 5.2.3, \cite{L5} Lemma 2.5, \cite{LY} Lemma 2.1]\label{characterization of deco tetra}
Suppose $\tau=\{1234\}$ is an ideal tetrahedron.
Then $(l_{12}, \cdots, l_{34})\in \mathbb{R}^6$ corresponds to the edge lengths of
a nondegenerate decorated hyperbolic ideal tetrahedron if and only if
\begin{equation}\label{triangle inequality}
\begin{aligned}
e^{(l_{ij}+l_{kh})/2}+e^{(l_{ik}+l_{jh})/2}>e^{(l_{ih}+l_{jk})/2}
\end{aligned}
\end{equation}
for $\{i,j,k,h\}=\{1,2,3,4\}$.
Furthermore, when $(l_{12}, \cdots, l_{34})\in \mathbb{R}^6$ is the edge length of a nondegenerate
hyperbolic ideal tetrahedron, the dihedral angle $\alpha_{ij}$ along the edge $v_iv_j$
equals the inner angle opposite to the edge with length $e^{(l_{ij}+l_{kh})/2}$
in the triangle formed by three edges with lengths $e^{(l_{ij}+l_{kh})/2}$, $e^{(l_{ik}+l_{jh})/2}$
and  $e^{(l_{ih}+l_{jk})/2}$.
\end{lemma}
Lemma \ref{characterization of deco tetra} implies that the dihedral angles along the edges in a decorated hyperbolic ideal tetrahedron
is independent of the choice of decorations and the dihedral angles along opposite edges are equal.
As a direct consequence of Lemma \ref{characterization of deco tetra},
the admissible space of edge lengths for a decorated hyperbolic ideal tetrahedron $\sigma$
is given by
\begin{equation*}
\begin{aligned}
\mathcal{L}_\sigma=\{(l_{12}, \cdots, l_{34})\in \mathbb{R}^6| &e^{(l_{ij}+l_{kh})/2}+e^{(l_{ik}+l_{jh})/2}>e^{(l_{ih}+l_{jk})/2},\\
                                               &\{i,j,k,h\}=\{1,2,3,4\}  \}.
\end{aligned}
\end{equation*}
The admissible space $\mathcal{L}_\sigma$ is nonconvex and invariant under the change of decorations.
Note that
$\mathbb{R}^6\setminus \mathcal{L}_\sigma=V_{12}^{34}\cup V_{13}^{24}\cup V_{14}^{23}$, where
\begin{equation*}
\begin{aligned}
V_{ij}^{kh}=&\{(l_{12}, \cdots, l_{34})\in \mathbb{R}^6| e^{(l_{ik}+l_{jh})/2}+e^{(l_{ih}+l_{jk})/2}\leq e^{(l_{ij}+l_{kh})/2}\}.
\end{aligned}
\end{equation*}
It is directly to see that $V_{ij}^{kh}=V_{kh}^{ij}$ and
\begin{equation*}
\begin{aligned}
V_{ij}^{kh}=&\{(l_{12}, \cdots, l_{34})\in \mathbb{R}^6| l_{ij}\geq -l_{kh}+2\ln [e^{(l_{ik}+l_{jh})/2}+e^{(l_{ih}+l_{jk})/2}]\}
\end{aligned}
\end{equation*}
is the closed region above an analytical function defined on $\mathbb{R}^5$.
Furthermore, $V_{12}^{34}, V_{13}^{24}, V_{14}^{23}$ are mutually disjoint, which follows from the fact that three positive constants $a,b,c$
can not satisfy $a\geq b+c$ and $b\geq a+c$ simultaneously. This implies the admissible space
$\mathcal{L}_\sigma=\mathbb{R}^6\setminus (V_{12}^{34}\cup V_{13}^{24}\cup V_{14}^{23})$
is homotopy equivalent to $\mathbb{R}^6$ and then simply connected.
In summary, we have the following corollary of Lemma \ref{characterization of deco tetra}.

\begin{corollary}\label{simply connectness}
The admissible space $\mathcal{L}_\sigma$ of edge lengths for a decorated hyperbolic ideal tetrahedron being nondegenerate is simply connected.
Furthermore, the boundary of $\mathcal{L}_\sigma$ consists of three disjoint analytical hypersurface of $\mathbb{R}^6$.
\end{corollary}

The idea to prove simply connectivity of admissible space by homotopy equivalence
was first introduced by the author in \cite{X1} and then further applied in \cite{HX1,HX2,X2}
to prove simply connectivity of admissible spaces in different cases.
It seems that this is a basic technique for handing similar problems.

For a decorated hyperbolic ideal tetrahedron $\sigma$, one has the following Schl\"{a}fli formula \cite{B1,L2a,Mil,R2,V},
\begin{equation}\label{Schlafli formula}
\begin{aligned}
-2d(vol_\sigma)=\sum_{i<j}l_{ij}d\alpha_{ij},
\end{aligned}
\end{equation}
where $vol_\sigma$ is the volume of the hyperbolic ideal tetrahedron.
The volume $vol_\sigma$ could be explicitly computed using the formula
$vol_\sigma=\frac{1}{2}\sum_{i<j}\Lambda(\alpha_{ij})$,
where $\Lambda(x)=-\int_0^x\ln |2\sin t|dt$ is the Lobachevsky function. One can refer to \cite{Mil,Rat} for this formula.

Dual to the volume, the co-volume of a decorated hyperbolic ideal tetrahedron $\sigma$ is defined to be
\begin{equation*}
\begin{aligned}
cov_\sigma=2vol_\sigma+\sum_{i<j}\alpha_{ij}l_{ij}.
\end{aligned}
\end{equation*}
By Schl\"{a}fli formula (\ref{Schlafli formula}), we have
\begin{equation*}
\begin{aligned}
d(cov_\sigma)=2d(vol_\sigma)+\sum_{i<j}d(\alpha_{ij}l_{ij})=\sum_{i<j}\alpha_{ij}dl_{ij},
\end{aligned}
\end{equation*}
which implies $\frac{\partial cov_\sigma}{\partial l_{ij}}=\alpha_{ij}$.
Then we have
\begin{equation*}
\begin{aligned}
\frac{\partial \alpha_{ij}}{\partial l_{kh}}=\frac{\partial^2 cov_\sigma}{\partial l_{kh} \partial l_{ij}}
=\frac{\partial^2 cov_\sigma}{\partial l_{ij} \partial l_{kh}}
=\frac{\partial \alpha_{kh}}{\partial l_{ij}},
\end{aligned}
\end{equation*}
which implies the Jacobian matrix $J_\sigma$ of the dihedral angles with respect to the edge lengths
in a decorated hyperbolic ideal tetrahedron is symmetric.
By Corollary \ref{simply connectness}, the following integral
\begin{equation}\label{F sigma}
\begin{aligned}
F_\sigma(l)=\int^l_{0}\sum_{i<j}\alpha_{ij}dl_{ij}
\end{aligned}
\end{equation}
is a well-defined smooth function on $\mathcal{L}_\sigma$.
As $dF_\sigma=\sum_{i<j}\alpha_{ij}dl_{ij}=d(cov_\sigma)$, $F_\sigma$ differs from the co-volume function $cov_\sigma$ by a constant.
Luo-Yang \cite{LY} found the deep relationship between the volume function $vol_\sigma$ and the co-volume function $cov_\sigma$ via Fenchel dual.
The readers can refer to \cite{LY} for more details on this.

For a decorated hyperbolic ideal tetrahedron
$\sigma=(\tau, \{H_1, H_2, H_3, H_4\})$,
using the relationship between edge lengths and inner angles in a Euclidean triangle $\tau\cap H_i$,
one has the following result on the Jacobian matrix $J_\sigma$ of the dihedral angles with respect
to the edge lengths.

\begin{lemma}[\cite{LY,Y}]\label{positivity of Jacobian for a tetra}
For a decorated hyperbolic ideal tetrahedron $\sigma$, the Jacobian matrix of the six dihedral angles
with respect to the six edge lengths is given by
\begin{equation*}
\begin{aligned}
J_\sigma=\left(\frac{\partial \alpha}{\partial l}\right)=\frac{1}{2}\left(
                                         \begin{array}{cc}
                                           M & M \\
                                           M & M \\
                                         \end{array}
                                       \right),
\end{aligned}
\end{equation*}
where $M$ is a $3\times 3$ matrix given by
\begin{equation*}
\begin{aligned}
M=\left(
    \begin{array}{ccc}
      \cot \alpha_{13}+\cot \alpha_{14} & -\cot \alpha_{14} & -\cot \alpha_{13} \\
      -\cot \alpha_{14} & \cot \alpha_{12}+\cot \alpha_{14} & -\cot \alpha_{12} \\
      -\cot \alpha_{13} & -\cot \alpha_{12} & \cot \alpha_{12}+\cot \alpha_{13} \\
    \end{array}
  \right).
\end{aligned}
\end{equation*}
Here the labels for the six columns and rows in $J_\sigma$
are $12, 13, 14, 34, 24, 23$ so that $i$-th column and $(i+3)$-th column, $i$-th row and $(i+3)$-th row
correspond to opposite edges.
Furthermore, $J_\sigma$ is positive semi-definite with null space spanned by the following $4$ linearly independent vectores
\begin{equation*}
\begin{aligned}
(1,1,1,0,0,0)^T, (1,0,0,0,1,1)^T, (0,1,0,1,0,1)^T,(0,0,1,1,1,0)^T.
\end{aligned}
\end{equation*}
\end{lemma}

Lemma \ref{positivity of Jacobian for a tetra} implies $F_\sigma$ is locally convex on $\mathcal{L}_\sigma$ and strictly locally convex on  $\mathcal{L}_\sigma\cap Ker(J_\sigma)^\perp$.

\subsection{Decorated hyperbolic polyhedral metrics on triangulated cusped 3-manifolds}
Suppose $(M, \mathcal{T})$ is an ideally triangulated cusped 3-manifold. For simplicity,
we denote the tetrahedra in $(M, \mathcal{T})$ as $\sigma_1,\cdots,\sigma_N$ in the following.

\begin{definition}
A decorated hyperbolic polyhedral metric on $(M, \mathcal{T})$ is a map $l: E\rightarrow \mathbb{R}$
such that for any ideal tetrahedron $\tau=\{ijkh\}\in \mathcal{T}$, the six real numbers $l_{ij}, l_{ik}, l_{ih}, l_{jk}, l_{jh}, l_{kh}$
form the edge lengths of a nondegenerate decorated hyperbolic ideal tetrahedron.
\end{definition}

We denote the admissible space of decorated hyperbolic polyhedral metrics on $(M, \mathcal{T})$
as $\mathcal{L}(\mathcal{T})$ in the following.
For simplicity of notations, we will label the edges with one index
when we are handling problems for ideally triangulated cusped 3-manifolds
and label the edges with two indices when we are handling problems for a single ideal tetrahedron in the following.

For the following applications, we introduce the following two notations.

\begin{definition}[\cite{FG,Y}]
The cusp relation matrix $C$ is an $s\times N$ matrix with the $(i,j)$-entry $c_{ij}$ equal to the number of ends of edge $j$ on the cusp $i$.
\end{definition}
The cusp relation matrix $C$ has the following properties.

\begin{proposition}[\cite{Choi,LY,N,Y}]\label{properties of C}
Suppose $C$ is the cusp relation matrix for an ideally triangulated cusped $3$-manifold $(M, \mathcal{T})$. Then
\begin{description}
  \item[(1)] $Rank(C)=s$, which implies $|V|\leq |E|=|T|$ for an  ideally triangulated cusped $3$-manifold $(M, \mathcal{T})$.
  \item[(2)] Two decorated hyperbolic polyhedral metrics $l_A, l_B$ correspond to the same hyperbolic polyhedral metric
             if and only if $l_A$ and  $l_B$ differ by a change of decorations, which is equivalent to $l_A-l_B\in Im(C^T)$.
\end{description}
\end{proposition}

\begin{definition}[\cite{Y}]
The incident matrix $G$ is the $N\times 6N$ matrix, so that for $n=1,\cdots, 6$ and $i,j=1,\cdots, N$, the $(i, 6j+n-6)$-entry of
$G$ is $1$ if and only if the $n$-th edge in the tetrahedron $\sigma_j$ is from the $i$-th edge in $E=E(\mathcal{T})$.
\end{definition}

Set $J=diag\{J_1, \cdots, J_N\}$, where $J_i$ is the Jacobian matrix for the tetrahedron $\sigma_i$.
Then $J$ is symmetric and positive semi-definite by Lemma \ref{positivity of Jacobian for a tetra}.
By the chain rules, the Jacobian matrix $\Lambda$ of the combinatorial Ricci curvature $K$
with respect to the edge length $l\in \mathcal{L}(\mathcal{T})$ is given by
\begin{equation}\label{expression of discrete Laplacian}
\begin{aligned}
\Lambda:=\left(\frac{\partial K}{\partial l}\right)=-GJG^T.
\end{aligned}
\end{equation}
The Jacobian matrix $\Lambda$ has the following useful property (see for instance \cite{Choi} Page 1354, \cite{Y} Proposition 4.13).

\begin{theorem}\label{negativity of L}
The Jacobian matrix $\Lambda$ of the combinatorial Ricci curvature $K$ with respect to the edge length $l\in \mathcal{L}(\mathcal{T})$
is symmetric and negative semi-definite.
Furthermore,
\begin{equation*}
\begin{aligned}
Ker(\Lambda)=Im(C^T)
\end{aligned}
\end{equation*}
and $\Lambda$ is strictly negative definite on $Ker(\Lambda)^\perp=Ker(C)$.
\end{theorem}

Theorem \ref{negativity of L} implies the decorated hyperbolic polyhedral metric
is locally determined by combinatorial Ricci curvature up to change of decorations.
The global rigidity of decorated hyperbolic polyhedral metric has been proved by Luo-Yang \cite{LY}.

Using the co-volume function $cov_{\sigma}$ for a single decorated hyperbolic ideal tetrahedron $\sigma$,
one can define the following co-volume function $cov$ for a
decorated hyperbolic polyhedral metric  $l\in \mathcal{L}(\mathcal{T})$
\begin{equation*}
\begin{aligned}
cov=\sum_{\sigma\in T}cov_\sigma,
\end{aligned}
\end{equation*}
which is a locally convex function defined on the admissible space $\mathcal{L}(\mathcal{T})$
with $\frac{\partial cov}{\partial l_i}=2\pi-K_i$
and strictly locally convex on $\mathcal{L}(\mathcal{T})\cap Ker(C)$.

\section{Basic properties of combinatorial Ricci flow}\label{Section 3}
Following Chow-Luo's combinatorial Ricci flow for circle packing metrics on surfaces \cite{CL}
and Luo's combinatorial Ricci flow for hyper-ideal polyhedral metrics on compact 3-manifolds with boundary \cite{L2},
Yang \cite{Y} introduced the following combinatorial
Ricci flow for decorated hyperbolic polyhedral metrics on cusped 3-manifolds.

\begin{definition}\cite{Y}
Suppose $(M, \mathcal{T})$ is an ideally triangulated cusped $3$-manifold.
The combinatorial Ricci flow for decorated hyperbolic polyhedral metrics on $(M, \mathcal{T})$ is defined to be
\begin{equation}\label{combinatorial Ricci flow context}
\begin{aligned}
\frac{dl_{i}}{dt}=K_{i}
\end{aligned}
\end{equation}
with $l(0)=l_0\in \mathcal{L}(\mathcal{T})$, where $l_{i}$ is the length of edge $i\in E$ and
$K_{i}$ is the combinatorial Ricci curvature along the edge $i$.
\end{definition}

The combinatorial Ricci flow (\ref{combinatorial Ricci flow context})
for decorated hyperbolic polyhedral metric is invariant under the change of decoration \cite{Y}.
In other words, if
$l(t)$ is a solution of the combinatorial Ricci flow (\ref{combinatorial Ricci flow context}) and
$\widetilde{l}(t)=l(t)+C^Tx$ for some $x\in \mathbb{R}^V$, then $\frac{d\widetilde{l}_i}{dt}=K_i(\widetilde{l})$,
which implies $\widetilde{l}(t)$ is also a solution of the combinatorial Ricci flow (\ref{combinatorial Ricci flow context}).
This follows from the fact that the dihedral angles and combinatorial Ricci curvatures are invariant under the change of decorations.
This also implies the combinatorial Ricci flow (\ref{combinatorial Ricci flow context})
is a well-defined deformation of the hyperbolic polyhedral metrics on $(M, \mathcal{T})$.
As the combinatorial Ricci flow (\ref{combinatorial Ricci flow context}) is an ODE system,
the short time existence of the solution of (\ref{combinatorial Ricci flow context}) follows
from the standard ODE theory \cite{Hart}.

\begin{lemma}[\cite{Y}]\label{invariance of sum of edge length}
Suppose $(M, \mathcal{T})$ is an ideally triangulated cusped $3$-manifold.
Then $\sum_{i\sim p}l_i$ is invariant along the combinatorial Ricci flow (\ref{combinatorial Ricci flow context}),
where $p\in V$ is a cusp of $(M, \mathcal{T})$ and $i\in E$ is an edge of $(M, \mathcal{T})$.
\end{lemma}
\proof
This follows from the fact that
\begin{equation*}
\begin{aligned}
\frac{d(\sum_{i\sim p}l_i)}{dt}=\sum_{i\sim p}K_i=0
\end{aligned}
\end{equation*}
for cusped 3-manifolds.
\qed

As a direct consequence of Lemma \ref{invariance of sum of edge length},
$\sum_{i=1}^Nl_i$ is invariant along the combinatorial Ricci flow (\ref{combinatorial Ricci flow context}).

\begin{remark}
By Lemma \ref{invariance of sum of edge length},
we can assume $\sum_{i\sim p}l_i(0)=0, \forall p\in V$, which is possible by Proposition \ref{properties of C}.
Then the solution of the combinatorial Ricci flow (\ref{combinatorial Ricci flow context}) stays in the plane
$\{l\in \mathbb{R}^N| \sum_{i\sim p}l_i=0, \forall p\in V \}$ by Lemma \ref{invariance of sum of edge length},
which is in fact $Ker(C)$.
For simplicity, we will assume $\sum_{i\sim p}l_i(0)=0$ in the following.
\end{remark}

Along the combinatorial Ricci flow (\ref{combinatorial Ricci flow context}),
the combinatorial Ricci curvature evolves according to the following equation
\begin{equation*}
\begin{aligned}
\frac{dK_i}{dt}=\sum_{j}\frac{\partial K_i}{\partial l_j}\frac{dl_j}{dt}=(\Lambda K)_i.
\end{aligned}
\end{equation*}
As $\Lambda$ is negative semi-definite, we can define the combinatorial Laplacian $\Delta$
for a decorated hyperbolic polyhedral metric $l$ as
$\Delta=\Lambda$. Then the combinatorial Ricci curvature evolves according to a discrete heat equation
$$\frac{dK}{dt}=\Delta K$$
along the combinatorial Ricci flow (\ref{combinatorial Ricci flow context}).
This is very similar to Hamilton's Ricci flow on 3-manifolds \cite{H}
and matches Luo's motivation to define a combinatorial Ricci flow for compact 3-manifolds with boundary \cite{L2}.

\begin{remark}
By Lemma \ref{positivity of Jacobian for a tetra} and the formula (\ref{expression of discrete Laplacian}),
the combinatorial Lapacian $\Delta$ is intrinsic in the sense
that it is independent of choice of decorations and depends only on the hyperbolic polyhedral metrics.
This implies the combinatorial Lapacian $\Delta$ is well-defined.
Using the combinatorial Laplacian $\Delta$, the author \cite{X4} introduced the following combinatorial Calabi flow
\begin{equation}\label{combinatorial Calabi flow}
\begin{aligned}
\frac{dl_i}{dt}=-\Delta K_i
\end{aligned}
\end{equation}
for decorated hyperbolic polyhedral metrics on cusped $3$-manifolds.
The basic properties of the combinatorial Calabi flow (\ref{combinatorial Calabi flow}), including
the stability of the combinatorial Calabi flow and others, were established in \cite{X4}.
For a single ideal tetrahedron, the author \cite{X4} further proved that
for any reasonable prescribed dihedral angles, the combinatorial Calabi flow (\ref{combinatorial Calabi flow})
exists for all time and converges exponentially fast to a decorated hyperbolic
polyhedral metric with the prescribed dihedral angles.
\end{remark}

The following property is a direct consequence of the smoothness of dihedral angles
in a decorated hyperbolic ideal tetrahedron with respect to the edge lengths.
\begin{lemma}[\cite{Y}]
If the solution of combinatorial Ricci flow (\ref{combinatorial Ricci flow context})
converges to some decorated hyperbolic polyhedral metric $l^*\in \mathcal{L}(\mathcal{T})$, then $K(l^*)=0$.
\end{lemma}

Furthermore, we have the following stability for the combinatorial Ricci flow (\ref{combinatorial Ricci flow context}) .
\begin{theorem}
Suppose $(M, \mathcal{T})$ is an ideally triangulated cusped $3$-manifold and
$l^*\in \mathcal{L}(\mathcal{T})\cap Ker(C)$ is a decorated hyperbolic polyhedral metric with $K(l^*)=0$.
Then there exists a positive constant $\delta$ such that if the initial value $l(0)=l_0\in \mathcal{L}(\mathcal{T})\cap Ker(C)$
of the combinatorial
Ricci flow (\ref{combinatorial Ricci flow context})  satisfies
$$||l_0-l^*||<\delta,$$
then the combinatorial Ricci flow (\ref{combinatorial Ricci flow context})  exists for all time and converges exponentially fast
to the decorated hyperbolic polyhedral metric $l^*$.
\end{theorem}
\proof
Set $\Gamma(l)=K(l)$. Take the combinatorial Ricci flow (\ref{combinatorial Ricci flow context}) as an autonomous system $\frac{dl}{dt}=\Gamma(l)$.
Then $l^*$ is a equilibrium point of the system by $K(l^*)=0$ and
\begin{equation*}
\begin{aligned}
D_l\Gamma|_{l^*}=\left(\frac{\partial K}{\partial l}\right)|_{l^*}=\Lambda(l^*)\leq 0
\end{aligned}
\end{equation*}
by Theorem \ref{negativity of L}.
Further note that the solution of the combinatorial Ricci flow (\ref{combinatorial Ricci flow context})
stays in $Ker(C)=Ker(\Lambda)^\perp$ by Lemma \ref{invariance of sum of edge length}
and $\Lambda$ is strictly negative definite on $Ker(C)=Ker(\Lambda)^\perp$ by Theorem \ref{negativity of L},
we have $l^*$ is a local attractor of the
combinatorial Ricci flow (\ref{combinatorial Ricci flow context}).
Then the conclusion follows from Lyapunov Stability Theorem (\cite{Pon}, Chapter 5).
\qed

For general initial values, the combinatorial Ricci flow may develop singularities, which corresponds to the boundary of
the admissible space $\mathcal{L}(\mathcal{T})$ \cite{Y}.
Under a very strong condition that the combinatorial Ricci flow does not develop singularities in finite and infinite time,
Yang \cite{Y} obtained the
the convergence of the combinatorial Ricci flow (\ref{combinatorial Ricci flow context}).
However, the condition that the combinatorial Ricci flow (\ref{combinatorial Ricci flow context})
does not develop singularities in finite and infinite time
is difficult to be satisfied in practice. We shall overcome this difficulty by extending
the combinatorial Ricci flow (\ref{combinatorial Ricci flow context}) through the potential singularities of the flow.

\section{Extension of the dihedral angles and the combinatorial Ricci curvature}\label{Section 4}
\subsection{Extension of dihedral angles}
Recall the admissible space of edge lengths for a decorated hyperbolic ideal tetrahedron $\sigma$ is given by
\begin{equation*}
\begin{aligned}
\mathcal{L}_\sigma=\mathbb{R}^6\setminus (V_{12}^{34}\cup V_{13}^{24}\cup V_{14}^{23}),
\end{aligned}
\end{equation*}
where each $V_{ij}^{kl}$, $\{i,j,k,l\}=\{1,2,3,4\}$, is
the closed region above an analytical function defined on $\mathbb{R}^5$ by Corollary \ref{simply connectness}.

\begin{lemma}\label{extension of dihedral angle}
Suppose $\sigma=(\tau, \{H_1, H_2, H_3, H_4\})$ is a decorated hyperbolic ideal tetrahedron
 and $\overline{l}\in \partial V_{ij}^{kh}$. If $l\in \mathcal{L}_\sigma$ tends to $\overline{l}$, then
$\alpha_{ij}=\alpha_{kh}\rightarrow \pi$ and $\alpha_{ik}, \alpha_{ih},\alpha_{jk},\alpha_{jh}\rightarrow 0$.
As a consequence, the dihedral angles of an ideal tetrahedron could be extended by constants
to be a continuous function defined  on $\mathbb{R}^6$.
\end{lemma}
\proof
For $l\in \mathcal{L}_\sigma$, by the cosine law for Euclidean triangle $\tau\cap H_i$, we have
\begin{equation*}
\begin{aligned}
\cos\alpha_{ij}
=\cos\alpha_{kh}
=\frac{y_{jk}^2+y_{jh}^2-y_{kh}^2}{2y_{jk}y_{jh}},
\end{aligned}
\end{equation*}
where $y_{jk}=e^{(l_{jk}+l_{ih})/2}, y_{jh}=e^{(l_{jh}+l_{ik})/2}, y_{kh}=e^{(l_{ij}+l_{kh})/2}$.
For $\overline{l}\in \partial V_{ij}^{kh}$,
we have $\overline{y}_{kh}=\overline{y}_{jk}+\overline{y}_{jh}$.
This implies $\cos\alpha_{ij}=\cos\alpha_{kh}\rightarrow -1$ as $l\rightarrow \overline{l}$.
Therefore, $\alpha_{ij}=\alpha_{kh}\rightarrow \pi$ as $l\rightarrow \overline{l}$.
$\alpha_{ik}, \alpha_{ih},\alpha_{jk},\alpha_{jh}\rightarrow 0$ can be proved similarly.

Define $\widetilde{\alpha}_{ij}$ as follows
\begin{equation*}
\begin{aligned}
\widetilde{\alpha}_{ij}=\left\{
                                                \begin{array}{ll}
                                                  \pi, & \hbox{$l\in V_{ij}^{kh}$;} \\
                                                  \alpha, & \hbox{$l\in \mathcal{L}_\sigma$;} \\
                                                  0, & \hbox{$l\in V_{ik}^{jh}\cup V_{ih}^{jk}$.}
                                                \end{array}
                                              \right.
\end{aligned}
\end{equation*}
Then $\widetilde{\alpha}_{ij}$ is a continuous function defined on $\mathbb{R}^6$ extending $\alpha_{ij}$.
The other dihedral angles can be extended similarly.
\qed

Note that for the extended dihedral angles, we still have $\widetilde{\alpha}_{ij}=\widetilde{\alpha}_{kh}$
and the sum of six extended dihedral angles is $2\pi$.
The extension of dihedral angles in Lemma \ref{extension of dihedral angle} is essentially Luo-Yang's extension in \cite{LY}.
Before going on, we recall the following definition and theorem of Luo in \cite{L4}.
\begin{definition}
A differential 1-form $w=\sum_{i=1}^n a_i(x)dx^i$ in an open set $U\subset \mathbb{R}^n$ is said to be continuous if
each $a_i(x)$ is continuous on $U$.  A continuous differential 1-form $w$ is called closed if $\int_{\partial \tau}w=0$ for each
triangle $\tau\subset U$.
\end{definition}

\begin{theorem}[\cite{L4}, Corollary 2.6]\label{Luo's convex extention}
Suppose $X\subset \mathbb{R}^n$ is an open convex set and $A\subset X$ is an open subset of $X$ bounded by a $C^1$
smooth codimension-1 submanifold in $X$. If $w=\sum_{i=1}^na_i(x)dx_i$ is a continuous closed 1-form on $A$ so that
$F(x)=\int_a^x w$ is locally convex on $A$ and each $a_i$ can be extended continuous to $X$ by constant functions to a
function $\widetilde{a}_i$ on $X$, then  $\widetilde{F}(x)=\int_a^x\sum_{i=1}^n\widetilde{a}_i(x)dx_i$ is a $C^1$-smooth
convex function on $X$ extending $F$.
\end{theorem}

By Theorem \ref{Luo's convex extention}, we have the following result.
\begin{lemma}\label{extension of F sigma}
The function $F_\sigma$ defined on $\mathcal{L}_\sigma$ in (\ref{F sigma})
could be extended to be a function
\begin{equation*}
\begin{aligned}
\widetilde{F}_\sigma(l)=\int_{0}^l\sum_{i<j}\widetilde{\alpha}_{ij}dl_{ij},
\end{aligned}
\end{equation*}
which is $C^1$ smooth and convex on $\mathbb{R}^6$ with $\nabla \widetilde{F}_\sigma=\widetilde{\alpha}$.
\end{lemma}
$\widetilde{F}_\sigma$ is the extended co-volume function up to a constant obtained by Luo-Yang \cite{LY} using Fenchel dual.
Here we give a direct extension of $F_\sigma$.
\begin{remark}
The idea of extension was first introduced by
Bobenko-Pinkall-Springborn \cite{BPS} to extend a locally convex function on a nonconvex domain to a convex function
and solved affirmably a conjecture of Luo \cite{L1} on the global rigidity of piecewise linear metrics on surfaces.
Luo \cite{L4} systematically studied the method of extension and
proved the global rigidity of inversive distance circle packing metrics for nonnegative inversive distance.
The method of extension has lots of applications, see \cite{GJ0,GJ1,GJ2,GJ3,GX3,GGLSW,GLSW,LY,X0,X1,X2} and others.
\end{remark}

\subsection{Extension of the combinatorial Ricci curvature}

Using the extension of dihedral angles in a decorated hyperbolic ideal tetrahedron,
we can extend the definition of combinatorial Ricci curvature for $l\in\mathbb{R}^N$ as follows
\begin{equation*}
\begin{aligned}
\widetilde{K}_{i}(l)=2\pi-\sum_{\sigma\in T} \widetilde{\alpha}_{i}(l),
\end{aligned}
\end{equation*}
where $\widetilde{\alpha}_i$ is the extension of the dihedral angle $\alpha_i$ along the edge $i\in E$ in a tetrahedron $\sigma$.
If $l\in \mathbb{R}^N$ corresponds to a nondegenerate decorated hyperbolic polyhedral metric in $\mathcal{L}(\mathcal{T})$,
then $\widetilde{K}=K$, which implies $\widetilde{K}$ is an extension of $K$.
Paralleling to the combinatorial Ricci curvature $K$, the extended combinatorial Ricci curvature $\widetilde{K}$ has the following property
on cusped $3$-manifolds.

\begin{lemma}\label{widetilde K in Ker C}
Suppose $(M, \mathcal{T})$ is an ideally triangulated cusped $3$-manifolds.
Then the extended combinatorial Ricci curvature $\widetilde{K}$ satisfies
$\sum_{i\sim p}\widetilde{K}_i=0$, $\forall p\in V$, which is equivalent to $\widetilde{K}\in Ker(C)$.
\end{lemma}

Following Luo-Yang \cite{LY}, we call $l\in \mathbb{R}^N$ as a generalized decorated hyperbolic polyhedral metric in the following.
Using the extension $\widetilde{F}_\sigma$ of $F_\sigma$, we
can define the following Ricci energy function
\begin{equation*}
\begin{aligned}
\widetilde{F}(l)=\sum_{\sigma\in T}\widetilde{F}_\sigma-2\pi\sum_{i=1}^Nl_i
\end{aligned}
\end{equation*}
for generalized decorated hyperbolic polyhedral metrics $l\in \mathbb{R}^N$.
By Lemma \ref{extension of F sigma},
we have
\begin{lemma}\label{C1 convexity of tilde F}
$\widetilde{F}(l)$ is a $C^1$ smooth convex function globally defined on $\mathbb{R}^N$ with $\nabla \widetilde{F}=-\widetilde{K}$.
\end{lemma}

\begin{remark}
One can also define $\widetilde{F}(l)$ directly by extending $-\int^l_0\sum_{i=1}^N K_idl_i$
to $-\int^l_0\sum_{i=1}^N\widetilde{K}_idl_i$, which is well-defined,
$C^1$ smooth and convex on $\mathbb{R}^N$ by Theorem \ref{Luo's convex extention} and Theorem \ref{negativity of L}.
\end{remark}

Using the Ricci energy function $\widetilde{F}$, one can prove the following rigidity for
generalized decorated hyperbolic polyhedral metric,
which has been obtained by Luo-Yang \cite{LY}.
However, our proof is direct, short and does not involve angle structures and Fenchel dual of the volume function.
For completeness, we present the proof here.

\begin{theorem}[\cite{LY}, Theorem 3.2]\label{rigidity of extended Ricci curvature}
Suppose a decorated hyperbolic polyhedral metric $l_A\in \mathcal{L}(\mathcal{T})$ and a generalized decorated hyperbolic polyhedral metric
$l_B\in \mathbb{R}^N$ have the same generalized combinatorial Ricci curvature $\widetilde{K}$.
Then $l_A$ and $l_B$ differ by a change of decoration.
\end{theorem}
\proof
Define
\begin{equation*}
\begin{aligned}
f(t)=\widetilde{F}(tl_B+(1-t)l_A)
\end{aligned}
\end{equation*}
for $t\in [0,1]$.
Then $f(t)$ is a $C^1$ smooth convex function on $[0,1]$ by the $C^1$ smoothness and convexity of $\widetilde{F}$ in Lemma \ref{C1 convexity of tilde F}.
By direct calculations, we have
\begin{equation*}
\begin{aligned}
f'(t)=&-\sum_{i=1}^N (l_{B,i}-l_{A,i})\widetilde{K}_i(tl_B+(1-t)l_A),
\end{aligned}
\end{equation*}
which gives
\begin{equation*}
\begin{aligned}
f'(0)=-\sum_{i=1}^N (l_{B,i}-l_{A,i})\widetilde{K}_i(l_A),\ \
f'(1)=-\sum_{i=1}^N (l_{B,i}-l_{A,i})\widetilde{K}_i(l_B).
\end{aligned}
\end{equation*}
By the condition that $l_A\in \mathcal{L}(\mathcal{T})$ and $l_B\in \mathbb{R}^N$ have the same generalized combinatorial Ricci curvature,
we have
$f'(0)=f'(1)$.
Note that $f$ is a $C^1$ smooth convex function, $f'(t)$ is an increasing function of $t\in [0,1]$.
Therefore, $f'(t)$ is a constant on $[0,1]$ by $f'(0)=f'(1)$.

Note that $\widetilde{K}$ is smooth around $l_A\in \mathcal{L}(\mathcal{T})$, which implies that $f(t)$
is smooth around $0$. Specially, $f(t)$ is a smooth function on $[0,\epsilon)$ for some $\epsilon\in (0,1)$.
Therefore, we have $f''(t)=0$ for $t\in [0,\epsilon)$.
By direct calculations, we have
\begin{equation*}
\begin{aligned}
f''(t)=&-(l_B-l_A)^T\Lambda(l_B-l_A)\equiv 0, t\in [0,\epsilon).
\end{aligned}
\end{equation*}
Then $l_B-l_A\in Ker(\Lambda)=Im(C^T)$ by Theorem \ref{negativity of L}.
Therefore, $l_B=l_A+C^Tx$ for some $x\in \mathbb{R}^k$,
which is equivalent to that $l_A$ and $l_B$ differ by a change of decoration.
\qed

\section{Extension of the combinatorial Ricci flow}\label{Section 5}
\subsection{Definition of the extended combinatorial Ricci flow}
Using the extension of combinatorial Ricci curvature,
we can define the extension of combinatorial Ricci flow for
generalized decorated hyperbolic polyhedral metric $l\in \mathbb{R}^N$  as
\begin{equation}\label{extended combinatorial Ricci flow}
\begin{aligned}
\frac{dl_{i}}{dt}=\widetilde{K}_{i}
\end{aligned}
\end{equation}
with $l(0)=l_0\in \mathbb{R}^N$.
\subsection{Uniqueness of the solution of extended combinatorial Ricci flow}
Although the term $\widetilde{K}$ on the righthand side of
the extended combinatorial Ricci flow (\ref{extended combinatorial Ricci flow}) is $C^0$ and not $C^1$ on $\mathbb{R}^N$,
we still have the uniqueness of the solution of
extended combinatorial Ricci flow (\ref{extended combinatorial Ricci flow}),
which is a consequence of the convexity of $C^1$ smooth Ricci
energy function $\widetilde{F}(l)$ defined on $\mathbb{R}^N$.

\begin{theorem}\label{uniqueness of extend CRF}
The solution of the extended combinatorial Ricci flow (\ref{extended combinatorial Ricci flow}) is unique for any initial generalized decorated
hyperbolic polyhedral metric $l(0)\in \mathbb{R}^N$.
\end{theorem}

\proof
The proof is paralleling to that of Theorem 5.1 in \cite{X3},
we only sketch the proof here. The idea of the proof comes from Ge-Hua \cite{GH}.

By Lemma \ref{C1 convexity of tilde F}, $\widetilde{F}$ is a $C^1$ smooth convex function on $\mathbb{R}^N$.
We can use the standard technique in PDE to mollify $\widetilde{F}$ to be $\widetilde{F}_\epsilon$, 
which is a $C^\infty$ convex function on $\mathbb{R}^N$
such that $\widetilde{F}_\epsilon \rightarrow \widetilde{F}$ in $C^1_{loc}(\mathbb{R}^N)$ as $\epsilon\rightarrow 0$.
For the smooth convex function $\widetilde{F}_\epsilon$, we have
\begin{equation}\label{uniqueness proof equ}
\begin{aligned}
\left(\nabla \widetilde{F}_\epsilon(l_1)-\nabla \widetilde{F}_\epsilon(l_2)\right)\cdot (l_1-l_2)\geq 0,
\end{aligned}
\end{equation}
which is a consequence of the monotonicity of the function 
$$f(t)=\left(\nabla \widetilde{F}_\epsilon(l_1)-\nabla \widetilde{F}_\epsilon(l_1+t(l_2-l_1))\right)\cdot (l_1-l_2).$$
Let $\epsilon\rightarrow 0$, (\ref{uniqueness proof equ}) gives
$\left(\nabla \widetilde{F}(l_1)-\nabla \widetilde{F}(l_2)\right)\cdot (l_1-l_2)\geq 0$,
which is equivalent to $\left(\widetilde{K}(l_1)-\widetilde{K}(l_2)\right)\cdot (l_1-l_2)\leq 0$.

If $l_1(t)$ and $l_{2}(t)$ are two solutions of the extended flow (\ref{extended combinatorial Ricci flow}) 
with the same initial value $l(0)\in \mathbb{R}^N$, then
for the function $g(t)=||l_1(t)-l_2(t)||^2$, we have $g(0)=0$, $g(t)\geq 0$ and
\begin{equation*}
\begin{aligned}
\frac{dg}{dt}=2(l_1-l_2)(\frac{dl_1}{dt}-\frac{dl_2}{dt})=2(l_1-l_2)(\widetilde{K}_1-\widetilde{K}_2)\leq 0,
\end{aligned}
\end{equation*}
which implies $g(t)\equiv 0$. Therefore, $l_1(t)\equiv l_2(t)$.
\qed

As a corollary of Theorem \ref{uniqueness of extend CRF},
we have the following result which shows that the solution of
the extended combinatorial Ricci flow (\ref{extended combinatorial Ricci flow}) extends the solution of the
combinatorial Ricci flow (\ref{combinatorial Ricci flow}) for any initial decorated hyperbolic polyhedral metric in $\mathcal{L}(\mathcal{T})$.

\begin{corollary}\label{context unique extension}
For any initial decorated hyperbolic polyhedral metric $l(0)\in \mathcal{L}(\mathcal{T})$,
denote the solutions of the combinatorial Ricci flow (\ref{combinatorial Ricci flow})
and the extended combinatorial Ricci flow (\ref{extended combinatorial Ricci flow})  as $l(t)$ and $\widetilde{l}(t)$ respectively.
Then
$\widetilde{l}(t)=l(t)$
whenever $l(t)$ exists in $\mathcal{L}(\mathcal{T})$.
\end{corollary}

\begin{remark}
Although the solution of the extended combinatorial Ricci flow (\ref{extended combinatorial Ricci flow})
is unique, there may exist some other different extensions of the solution of the combinatorial Ricci flow (\ref{combinatorial Ricci flow}) on $\mathbb{R}^N$.
The uniqueness of extension of the solution of combinatorial Ricci flow  (\ref{combinatorial Ricci flow})
depends on the choice of extension of the combinatorial Ricci curvature.
Here we choose a natural extension of the combinatorial Ricci curvature, but there may exist some other different extension
of the combinatorial Ricci curvature on $\mathbb{R}^N$.
The author thanks Tian Yang
for pointing this out to the author.
\end{remark}

\subsection{Longtime existence of the extended combinatorial Ricci flow}
We have the following longtime existence for the solution of the extended combinatorial Ricci flow (\ref{extended combinatorial Ricci flow}).

\begin{theorem}\label{long time existence}
The solution of the extended combinatorial Ricci flow (\ref{extended combinatorial Ricci flow}) exists for all time
for any initial generalized decorated
hyperbolic polyhedral metric $l(0)\in \mathbb{R}^N$.
\end{theorem}
\proof
To prove the longtime existence of the solution of the extended combinatorial Ricci flow (\ref{extended combinatorial Ricci flow}),
we just need to prove that the solution $l(t)$ is bounded for any finite time $[0, T]$ with $T<+\infty$.

As the extended dihedral angles are uniformly bounded by $\pi$ and we are working 
on a cusped $3$-manifold $M$ with a fixed triangulation $\mathcal{T}$,
the extended combinatorial Ricci curvature $\widetilde{K}$ is uniformly bounded by some constant $C$
depending only on the triangulation $\mathcal{T}$, which implies that the solution $l(t)$ of
the extended combinatorial Ricci flow (\ref{extended combinatorial Ricci flow})
satisfies
\begin{equation*}
\begin{aligned}
|l_i(t)|\leq |l_i(0)|+CT<+\infty
\end{aligned}
\end{equation*}
for any $t\in [0, T]$.
Therefore, the solution of the extended combinatorial Ricci flow (\ref{extended combinatorial Ricci flow})
exists for all time for any initial generalized decorated
hyperbolic polyhedral metric $l(0)\in \mathbb{R}^N$.
\qed

\subsection{Convergence of the extended combinatorial Ricci flow}
\begin{lemma}\label{limit of extended CRF}
If the solution $l(t)$ of the extended combinatorial Ricci flow (\ref{extended combinatorial Ricci flow})
on an ideally triangulated cusped $3$-manifold $(M, \mathcal{T})$
converges to some generalized decorated hyperbolic polyhedral metric $\overline{l}\in \mathbb{R}^N$, then $\widetilde{K}(\overline{l})=0$.
\end{lemma}
\proof
As $l(t)$ is a solution of (\ref{extended combinatorial Ricci flow}), which exists for all time by Theorem \ref{long time existence},
there exists $\xi_n\in (n, n+1)$ such that
\begin{equation*}
\begin{aligned}
l(n+1)-l(n)=l'(\xi_n)=\widetilde{K}(l(\xi_n))\rightarrow 0, n\rightarrow \infty,
\end{aligned}
\end{equation*}
by the convergence of the solution $l(t)$.
As $l(t)$ converges to $\overline{l}$, we have $l(\xi_n)\rightarrow \overline{l}$.
By the continuity of $\widetilde{K}$, we have $\widetilde{K}(\overline{l})=0$.
\qed

The following result is a generalization of Theorem \ref{main theorem introduction},
which gives a new characterization of the existence of a complete hyperbolic metric on an ideally triangulated cusped 3-manifold $(M, \mathcal{T})$
dual to Casson and Rivin's program.

\begin{theorem}\label{theorem equivalence sec 5}
Suppose there exists a complete hyperbolic metric on the ideally triangulated cusped $3$-manifold $(M, \mathcal{T})$.
Then $l^*\in \mathcal{L}(\mathcal{T})$ is a decorated hyperbolic polyhedral metric with zero combinatorial Ricci curvature
if and only if the solution of
extended combinatorial Ricci flow (\ref{extended combinatorial Ricci flow})
with initial value $l(0)\in \mathbb{R}^N$ exists for all time and converges exponentially fast to $l^*$.
\end{theorem}

\proof
Suppose the solution of extended combinatorial Ricci flow (\ref{extended combinatorial Ricci flow})
exists for all time and converges exponentially fast to $l^*$, then $\widetilde{K}(l^*)=0$ by Lemma \ref{limit of extended CRF}.
As there exists a nondegenerate decorated hyperbolic polyhedral metric with zero combinatorial Ricci curvature
on the ideal triangulated cusped 3-manifold $(M, \mathcal{T})$ by assumption, we have $l^*$ is
the nondegenerate decorated hyperbolic polyhedral metric with zero combinatorial Ricci curvature by Theorem \ref{rigidity of extended Ricci curvature}.
Therefore, $l^*\in \mathcal{L}(\mathcal{T})$ and $K(l^*)=0$.

Suppose $l^*\in \mathcal{L}(\mathcal{T})$ is a decorated hyperbolic polyhedral metric with zero combinatorial Ricci curvature.
Along the extended combinatorial Ricci flow (\ref{extended combinatorial Ricci flow}), we still have
\begin{equation*}
\begin{aligned}
\frac{d(\sum_{i\sim p}l_i)}{dt}=\sum_{i\sim p} \widetilde{K}_i=0
\end{aligned}
\end{equation*}
by Lemma \ref{widetilde K in Ker C},
which implies $\sum_{i\sim p}l_i$ is invariant along
the extended combinatorial Ricci flow (\ref{extended combinatorial Ricci flow}) for any cusp $p\in V$.
Without loss of generality, we assume that the initial decorated hyperbolic polyhedral metric $l(0)$ and $l^*$ are in $Ker(C)$,
where
\begin{equation*}
\begin{aligned}
Ker(C)=\{l\in \mathbb{R}^N|\sum_{i\sim p}l_i=0, \forall p\in V\}.
\end{aligned}
\end{equation*}
Then the solution $l(t)$ of the extended combinatorial Ricci flow (\ref{extended combinatorial Ricci flow}) stays in $Ker(C)$.

If $l^*\in \mathcal{L}(\mathcal{T})$ is a decorated hyperbolic polyhedral metric with $K(l^*)=0$,
then we have $\nabla \widetilde{F}(l^*)=-\widetilde{K}(l^*)=0$.
As $\widetilde{F}(l)$ is a $C^1$ smooth convex function on $\mathbb{R}^N$,
we have $l^*$ is a minimal point of $\widetilde{F}$ on $\mathbb{R}^N$.
Further note that $Hess \widetilde{F}(l^*)=-\Lambda(l^*)$ is positive semi-definite with
$Ker(Hess \widetilde{F}(l^*))^\perp=Ker(\Lambda)^\perp =Ker (C)$ by Theorem \ref{negativity of L},
$\widetilde{F}$ is a convex function on $Ker(C)$ and strictly locally convex around $l^*$ in $Ker(C)\cap \mathcal{L}(\mathcal{T})$.
Then we have $\lim_{l\in Ker(C),l\rightarrow \infty}\widetilde{F}(l)=+\infty$,
which implies $\widetilde{F}|_{Ker(C)}$ is proper.
This follows from the following result of convex function, the proof of which could be found in \cite{GX3} (Lemma 4.6).

\begin{lemma}
Suppose $f(x)$ is a $C^1$ smooth convex function on $\mathbb{R}^n$ with $\nabla f(x_0)=0$ for some $x_0\in \mathbb{R}^n$,
$f(x)$ is $C^2$ smooth and strictly convex
in a neighborhood of $x_0$, then $\lim_{x\rightarrow \infty}f(x)=+\infty$.
\end{lemma}

For any initial value $l(0)\in \mathbb{R}^N$, we have
\begin{equation}\label{derivative of tilde F}
\begin{aligned}
\frac{d\widetilde{F}(l(t))}{dt}=-\sum_i \widetilde{K}^2_i\leq 0,
\end{aligned}
\end{equation}
which implies $\widetilde{F}(l(t))$ is bounded along the extended combinatorial Ricci flow (\ref{extended combinatorial Ricci flow}).
Combining with the properness of convex function $\widetilde{F}|_{Ker(C)}$,
the solution $l(t)$ of the extended combinatorial Ricci flow (\ref{extended combinatorial Ricci flow})
lies in a compact subset $\Omega\subset\subset Ker(C)$.

By the fact that $\widetilde{F}(l(t))$ is bounded and (\ref{derivative of tilde F}), the limit $\lim_{t\rightarrow +\infty}\widetilde{F}(l(t))$ exists.
This implies
\begin{equation*}
\begin{aligned}
\widetilde{F}(l(n+1))-\widetilde{F}(l(n))=-\sum_i \widetilde{K}^2_i(l(\xi_n))\rightarrow 0, n\rightarrow \infty,
\end{aligned}
\end{equation*}
for some $\xi_n\in (n, n+1)$.
By the compactness of $\Omega$, we have a subsequence of $l(\xi_n)$ converges to some $\xi^*\in \Omega\subset Ker(C)$.
The continuity of $\widetilde{K}$ then implies $\widetilde{K}(\xi^*)=0$.
Therefore, we have
$$\widetilde{K}(l^*)=\widetilde{K}(\xi^*)=0$$
with $l^*\in \mathcal{L}(\mathcal{T})\cap Ker C$ and $\xi^*\in Ker C$.
By the rigidity of the extended combinatorial Ricci curvature $\widetilde{K}$ in Theorem \ref{rigidity of extended Ricci curvature}, we have $l^*=\xi^*$.
Therefore, $l(\xi_n)\rightarrow l^*$.

At $l^*$, we have
\begin{equation*}
\begin{aligned}
\frac{\partial \widetilde{K}}{\partial l}|_{l=l^*}=\frac{\partial K}{\partial l}|_{l=l^*}=\Lambda(l^*)\leq 0.
\end{aligned}
\end{equation*}
Furthermore, the solution $l(t)$ stays in $Ker(C)=Ker(\Lambda)^\perp$
along the extended combinatorial Ricci flow (\ref{extended combinatorial Ricci flow})
and $\Lambda(l^*)$ is strictly negative definite on $Ker(C)$ by Theorem \ref{negativity of L}.
This implies $l^*$ is a local attractor of the extended combinatorial Ricci flow (\ref{extended combinatorial Ricci flow}).
As $l(\xi_n)\rightarrow l^*$, it follows from Lyapunov Stability Theorem (\cite{Pon} Chapter 5) that
the solution $l(t)$ of the extended combinatorial Ricci flow (\ref{extended combinatorial Ricci flow}) converges exponentially fast to $l^*$.
\qed

\begin{remark}
By the proof of Theorem \ref{theorem equivalence sec 5}, $\sum_{i\sim p}l_i$ is invariant along
the extended combinatorial Ricci flow (\ref{extended combinatorial Ricci flow}) for any cusp $p\in V$,
which implies $\sum_{i=1}^N l_i$ is a constant along the extended combinatorial Ricci flow (\ref{extended combinatorial Ricci flow}).
Therefore, the Ricci energy $\widetilde{F}$ differs from the co-volume $cov$ by a constant along
the extended combinatorial Ricci flow (\ref{extended combinatorial Ricci flow}).
Note that the extended combinatorial Ricci flow (\ref{extended combinatorial Ricci flow}) finds the complete
hyperbolic metric on the cusped $3$-manifold $(M, \mathcal{T})$ by minimizing the Ricci energy $\widetilde{F}$,
this implies that the extended combinatorial Ricci flow (\ref{extended combinatorial Ricci flow}) 
find the complete hyperbolic metric on $(M, \mathcal{T})$ also by minimizing the co-volume function $cov$
along the flow.
This is dual to Casson-Rivin's program to find complete hyperbolic metrics on $(M, \mathcal{T})$
by maximizing the volume on angle structures.
\end{remark}

\textbf{Acknowledgements}\\[8pt]
Part of the series of work, including \cite{X3,X4} and this paper, was done when the author was visiting the Rutgers University.
The author thanks Professor Feng Luo for his invitation to Rutgers University and communications and interesting on this work.
The author thanks Professor Feng Luo for explaining the details of his joint work \cite{LY} with Tian Yang
and introducing the reference \cite{Y} to the author.
Part of the series of work was reported in an online seminar of Rutgers University in April 2020.
The author thanks participants in the seminar for communications and suggestions.
The research of the author is supported by
Fundamental Research Funds for the Central Universities
and National Natural Science Foundation of China under grant no. 61772379.

(Xu Xu) School of Mathematics and Statistics, Wuhan University, Wuhan 430072, P.R. China

E-mail: xuxu2@whu.edu.cn\\[2pt]

\end{document}